\documentclass{amsart}

\usepackage{amsmath}
\usepackage{amscd}
\usepackage{amssymb}

\newcommand{\bC}{{\mathbb C}}

\newcommand{\bP}{{\mathbb P}}

\newcommand{\bZ}{{\mathbb Z}}

\newcommand{\cE}{{\mathcal E}}
\newcommand{\cM}{{\mathcal M}}
\newcommand{\cO}{{\mathcal O}}
\newcommand{\cP}{{\mathcal P}}

\newcommand{\cW}{{\mathcal W}}
\newcommand{\Mbar}{\overline{\cM}}
\newcommand{\vac}{|0\rangle}
\newcommand{\lvac}{\langle 0|}

\newcommand{\MgX}{\Mbar_{g,0}(X,\mu^+,\mu^-)}
\newcommand{\gmu}[1]{\mathcal{#1}_{g,\mu^+,\mu^-}}

\newtheorem{conjecture}{Conjecture}

\newtheorem{theorem}{Theorem}[section]
\newtheorem{theorem/definition}{Theorem/Definition}[section]

\newtheorem{proposition}{Proposition}[section]
\newtheorem{lemma}{Lemma}[section]

\theoremstyle{remark}
\newtheorem{remark}{Remark}[section]

\theoremstyle{definition}

\begin{document}

\title
{Hodge Integrals and Integrable Hierarchies}
\author{Jian Zhou}
\address{Department of Mathematical Sciences\\Tsinghua University\\Beijing, 100084, China}
\email{jzhou@math.tsinghua.edu.cn}
\begin{abstract}
We show that the generating series of some Hodge integrals involving
one or two partitions are $\tau$-functions
of the KP hierarchy or the $2$-Toda hierarchy respectively.
We also formulate a conjecture on the connection between relative invariants and integrable hierarchies.
The conjecture is verified in some examples.
\end{abstract}
\maketitle

\section{Introduction}

Recently there have been some progresses in the study of Hodge integrals.
These are integrals of the form
$$\int_{\Mbar_{g, n}} \psi_1^{j_1}
\cdots \psi_n^{j_n}\lambda_1^{k_1} \cdots \lambda_g^{k_g}$$
where
$\psi_i$ and $\lambda_j$ are Chern classes of some bundles naturally defined on
the Deligne-Mumford space $\Mbar_{g, n}$.
A famous conjecture of Witten \cite{Wit} proved by Kontsevich \cite{Kon}
states that the generating series of integrals of the form
$$\int_{\Mbar_{g, n}} \psi_1^{j_1} \cdots \psi_n^{j_n}$$
is a $\tau$-function of the KdV hierarchy.
This yields recursion relations that completely determine all the integrals of this form.
As discovered in \cite{Zho1} and illustrated in the proofs of some formulas for some special
Hodge integrals that involve one or two partitions \cite{LLZ1, LLZ2, LLZ3},
a single partial differential equation called the {\em cut-and-join equation}
is a very effective tool to study such Hodge integrals.

There have been some interests in the relationship
between the cut-and-join equation and integrable hierarchies.
The purpose of this paper is to establish a relationship between the one-partition Hodge integrals
in \cite{Kat-Liu, Mar-Vaf, Zho1, LLZ1, LLZ2} and KP hierarchy,
and a relationship between the two-partition Hodge integrals in \cite{Zho2, LLZ3} and
the $2$-Toda hierarchy.
Our results are based on some formulas expressing Hodge integrals in terms of
some combinatorial objects,
such as Chern-Simons link invariants or modular matrices in the representation theory of Kac-Moody algebras
and conformal field theory.
In the case of one-partition Hodge integrals,
the formula was conjectured by Mari\~no and Vafa \cite{Mar-Vaf} and proved in joint work with
C.-C. Liu and K. Liu \cite{LLZ1, LLZ2};
in the case of two-partition Hodge integrals,
the formula was conjectured by the author \cite{Zho2} and proved in joint work also with
C.-C. Liu and K. Liu \cite{LLZ3}.
By the combinatorial studies in \cite{Zho1} and \cite{Zho2} respectively,
it is possible to express one-partition and two-partition Hodge integrals in terms of
Schur and skew Schur functions respectively
(see also Section \ref{sec:W} and Section \ref{sec:W2}).
Following the approach developed by the Kyoto school \cite{MJD},
it is then possible to rewrite the Hodge integrals as vacuum expectation values of some
operators on the fermionic Fock space,
and so they are $\tau$-functions of the corresponding integrable hierarchies by standard arguments.

Our work leaves still open the relationship between the cut-and-join equation
with integrable hierarchies.
Earlier result \cite{Fre-Wan} has revealed the relationship between the cut-and-join operator
and the Virasoro algebra,
but not to integrable hierarchies.
It might be the case that the cut-and-join equation are related to the Virasoro constraint.
We leave this for future investigations.
Finally it is interesting problem
to reprove Witten's conjecture/Kontsevich's theorm from the results in this work.

After the completion of this paper, Vafa informs the author that he has also
obtained the results in this work and results for the topological vertex.

The rest of the paper is arranged as follows.
In Section \ref{sec:Pre} we recall some well-known facts about integrable hierarchies.
In Section \ref{sec:One} we show the generating series of one-partition Hodge integrals
is a $\tau$-function of the KP hierarchy.
In Section \ref{sec:Two} we show the generating series of two-partition Hodge integrals
gives rise to a sequence of $\tau$-functions of the $2$-Toda hierarchy.
In Section \ref{sec:Relative} we extend our results to relative invariants for some examples.
A conjecture in the general case is formulated in this section.

\section{Preliminaries on Integrable Hierarchies}
\label{sec:Pre}

In this section we recall some well-known facts about integrable hierarchies.
For more details, see e.g. \cite{MJD, Kac, UT}.

\subsection{The charged free fermions}

The {\em charged free fermions} is a Lie superalgebra spanned by odd generators
$\{\psi^{\pm}_r: r\in \frac{1}{2}+\bZ\}$ and an even central generator $I$,
with the following commutation relations:
\begin{align} \label{eqn:CommF}
[\psi^+_r, \psi^-_s] & = \delta_{r, -s}I, &
[\psi^{\pm}_r, \psi^{\pm}_s] & = 0.
\end{align}
Consider the fermionic Fock space $F$ spanned by elements of the form:
\begin{eqnarray} \label{eqn:Basis}
&& \{\psi^+_{r_1} \cdots \psi^+_{r_m} \psi^-_{s_1} \cdots \psi^-_{s_n}| \;
m, n \geq 0, r_1 < \dots < r_m < 0, s_1 <  \dots < s_n < 0\}.
\end{eqnarray}
This space contains an element $1$,
which corresponds to the case $m=n=0$.
We denote this vector by $\vac$.
For $n \in \bZ -\{0\}$,
define
\begin{eqnarray}
|n\rangle & = & \begin{cases}
\psi^+_{-n+\frac{1}{2}} \cdots \psi^+_{-\frac{3}{2}}\psi^+_{-\frac{1}{2}}, & n > 0, \\
\psi^-_{n+\frac{1}{2}} \cdots \psi^-_{-\frac{3}{2}}\psi^-_{-\frac{1}{2}}, & n < 0.
\end{cases}
\end{eqnarray}

Let $I$ act on $F$ as the identity,
and $\psi^{\pm}_r$ act on $F$ by multiplication by $\psi^{\pm}_r$ when $r < 0$,
and by contraction with $\psi^{\mp}_{-r}$ when $r > 0$.
Then the commutation relations (\ref{eqn:CommF}) hold.
Define fields:
\begin{align*}
\psi^{\pm}(z) = & \sum_{r \in \frac{1}{2} + \bZ} \psi^{\pm}_r z^{-r - \frac{1}{2}}.
\end{align*}
Then one has the following OPE:
\begin{eqnarray}
&& \psi^+(z)\psi^-w) \sim \frac{1}{z-w}, \\
&& \psi^{\pm}(z)\psi^{\pm}(w) \sim 0.
\end{eqnarray}

One can define a Hermitian metric on $F$ by taking (\ref{eqn:Basis}) as orthonormal basis.
The {\em vacuum expectation value} (vev) of an operator $A: F \to F$ is defined by:
$$\langle A\rangle = \lvac A \vac.$$
Here we have followed the physicists' notation in writing the inner product of a vector $|v\rangle$
with a vector $|w\rangle$ as $\langle w|v\rangle$.

\subsection{Charge decomposition of $F$}

Consider the field:
$$\alpha(z) =:\psi^+(z)\psi^-(z):.$$
Then one has the following OPE:
\begin{eqnarray}
&& \alpha(z)\psi^{\pm}w) \sim \pm \frac{\psi^{\pm}(w)}{z-w}, \\
&& \alpha(z)\alpha(w) \sim \frac{1}{(z-w)^2}. \label{eqn:OPEalpha/alpha}
\end{eqnarray}
Write
$$\alpha(z) = \sum_{n \in \bZ} \alpha_n z^{-n-1}.$$
Then (\ref{eqn:OPEalpha/alpha}) is equivalent to:
\begin{eqnarray}
&& [\alpha_m, \alpha_n] + m\delta_{m, -n}.
\end{eqnarray}
The operator
$$\alpha_0 = \sum_{r \in \frac{1}{2}} :\psi^+_r\psi^-{_r}:$$
is called the {\em charge operator}.
It is diagonalizable on $F$:
the vector
$$\psi^+_{r_1} \cdots \psi^+_{r_m} \psi^-_{s_1} \cdots \psi^-_{s_n}$$
has eigenvalue $m-n$.
One gets a decomposition:
$$F = \oplus_{n \in \bZ} F^{(n)},$$
where $F^{(n)}$ is the eigenspace of $\alpha_0$ on which $\alpha_0$ has eigenvalue $n$.
An operator $A: F \to F$ is said to have charge $0$ if it preserves this decomposition.

Another way to describe the fermionic Fock space is by the semi-infinite wedge
(cf. e.g. \cite{Kac-Rai}).
It is very easy in this description to see that every partition $\mu$
corresponds to a vector $|\mu\rangle$ in $F^{(0)}$.
In particular the empty partition $(0)$ corresponds to the vacuum vector $\vac$.
The vectors $\{|\mu\rangle: \mu \in \cP\}$ form a basis of $F^{(0)}$.
Here $\cP$ denotes the set of all partitions, including $(0)$.

\subsection{The $\tau$-functions of the KP hierarchy as vevs}

For $t=(t_1, t_2, \dots)$,
let
$$\Gamma_{\pm}(t) = \exp(\sum_{n \geq 1} t_n \alpha_{\pm n}).$$

\begin{proposition} \label{prop:KP}
Suppose $A: F \to F$ is an operator of charge $0$,
which satisfies:
\begin{eqnarray} \label{eqn:Hirota1}
&& [A \otimes A, \sum_{r \in \frac{1}{2}+\bZ}\psi^+_r \otimes \psi^-_{-r}] = 0,
\end{eqnarray}
then
$$\tau(t) = \lvac\Gamma_+(t)A\vac$$
is a $\tau$-function of the KP hierarchy.
\end{proposition}

\begin{proof}
Recall if $|v\rangle \in F^{(0)}$ satisfies
\begin{eqnarray} \label{eqn:Hirota2}
\sum_{r\in \frac{1}{2} + \bZ} \psi^+_r|v\rangle \otimes \psi^-_{-r}|v\rangle = 0,
\end{eqnarray}
then
$$\tau(t) = \lvac\Gamma_+(t)|v\rangle$$ is a $\tau$-function of the KP hierarchy.
For a proof see e.g. \cite{MJD}.
Now if $A$ satisfies (\ref{eqn:Hirota1}),
then it is straightforward to see that
$A\vac$ satisfies (\ref{eqn:Hirota2}).
\end{proof}

\begin{remark}
If $A$ lies in $GL(\infty)$ then (\ref{eqn:Hirota1}) is automatically satisfied.
\end{remark}

\subsection{The $\tau$-functions of the $2$-Toda hierarchy as vevs}

Similar to Proposition \ref{prop:KP} one has (cf. \cite{UT} and \cite[Appendix]{Oko}):

\begin{proposition} \label{prop:Toda}
Suppose $A: F \to F$ is an operator of charge $0$,
which satisfies (\ref{eqn:Hirota1}),
then
\begin{eqnarray}
&& \tau_n(t^+, t^-) = \langle n| \Gamma_+(t^+)A \Gamma_-(t^-) |n\rangle
\end{eqnarray}
is a sequence of $\tau$-functions for the $2$-Toda hierarchy.
\end{proposition}

\subsection{Schur functions and skew Schur functions as vacuum expectation values}
For definitions of Schur functions and skew Schur functions see e.g. \cite{Mac}.

For formal variables $x=(x_1, x_2, \dots)$,
let
\begin{eqnarray*}
&& Y_{\pm}(x_1, x_2, \dots, )
= \exp\left( \sum_{n > 0} \frac{p_n(x_1, x_2, \dots)}{n} \alpha_{\pm n}\right),
\end{eqnarray*}
where
$$p_n(x_1, x_2, \dots) = \sum_i x_i^n$$
is the $n$-th Newton function.
Then clearly one has
\begin{eqnarray} \label{eqn:YGamma}
&& Y_{\pm}(x) = \Gamma_{\pm}(t),
\end{eqnarray}
for
$$t_n = \frac{1}{n}p_n(x).$$
One has (cf. e.g. \cite[Appendix]{Oko}):
\begin{eqnarray}
&& \langle\mu|Y_-(x)|\nu\rangle = \langle\nu|Y_+(x)|\mu\rangle =
s_{\mu/\nu}(x). \label{eqn:Skew}
\end{eqnarray}
In particular,
\begin{eqnarray}
&& \langle\mu|Y_-(x)|0\rangle = \langle 0|Y_+(x)|\mu\rangle =
s_{\mu}(x). \label{eqn:Schur}
\end{eqnarray}

For a partition $\mu = (\mu_1 \geq \mu_2 \geq \cdots \geq \mu_h > 0)$,
define
$$\kappa_{\mu} = \sum_{i=1}^h \mu_i(\mu_i-2i+1).$$
For the empty partition $(0)$,
we define
$$\kappa_{(0)} = 0.$$
Note we have
\begin{eqnarray*}
&& \kappa_{\mu} = \sum_i [(\mu_i-i+\frac{1}{2})^2- (-i+\frac{1}{2})^2].
\end{eqnarray*}

Finally, the operator
$$K=\frac{1}{2} \sum_{r \in \frac{1}{2}+\bZ} r^2:\psi^+_r\psi_{-r}^-:$$
has the following action on base vectors \cite{Oko1}:
\begin{eqnarray} \label{eqn:K}
K |\nu\rangle = \frac{\kappa_{\nu}}{2}|\nu\rangle;
\end{eqnarray}
furthermore,
let $R$ be the translation operator,
then
\begin{eqnarray} \label{eqn:Translation}
&& R^{-n}KR^n = K + n H +\frac{n^2}{2}\alpha_0 + \frac{n(4n^2-1)}{24},
\end{eqnarray}
where
$$H = \sum_{r \in \frac{1}{2}+\bZ} r :\psi^+_r\psi_{-r}^-: $$
is the energy operator,
whose action on the basis $\{|\mu\rangle\}$ is given by:
\begin{eqnarray} \label{eqn:H}
H|\mu\rangle = |\mu| \;|\mu\rangle.
\end{eqnarray}
The following commutation relation is also well-known:
\begin{eqnarray} \label{eqn:Comm}
&& [R, \Gamma_{\pm}(t)] = 0.
\end{eqnarray}

\section{One-Partition Hodge Integrals and KP Hierarchy}

\label{sec:One}

\subsection{One-partition Hodge Integrals}

For a partition $\mu$,
consider the following generating series of Hodge integrals:
\begin{eqnarray*}
G_{\mu}(r; \lambda)
& = & - \frac{\sqrt{-1}^{l(\mu)}}{z_{\mu^+}}
\cdot \left[r(r+1)\right]^{l(\mu)-1}
\cdot \prod_{i=1}^{l(\mu)} \frac{\prod_{a=1}^{\mu_i-1}
\left( \mu_ir + a \right)}{\mu_i!}\\
&& \cdot \sum_{g \geq 0} \lambda^{2g-2} \int_{\Mbar_{g, l(\mu)}}
\frac{\Lambda_{g}^{\vee}(1)\Lambda^{\vee}_{g}(r)\Lambda_{g}^{\vee}(-1 - r)}
{\prod_{i=1}^{l(\mu)} \frac{1}{\mu_i} \left(\frac{1}{\mu_i} - \psi_i\right)}.
\end{eqnarray*}
The following formula was conjectured by Mari\~no and Vafa \cite{Mar-Vaf} and proved in joint work
with C.-C. Liu and K. Liu \cite{LLZ1, LLZ2}:
\begin{eqnarray} \label{eqn:MV}
&& G^{\bullet}(r; \lambda; p) = R^{\bullet}(r;p),
\end{eqnarray}
where
\begin{eqnarray*}
G^{\bullet}(r; \lambda; p) & = &\exp \left(
\sum_{(\mu) \in \cP_+} G_{\mu}(r; \lambda)p_{\mu}\right), \\
R^{\bullet}(r; \lambda; p) & = & \sum_{\mu, \nu} \frac{\chi_{\nu}(\mu)}{z_{\mu}}
e^{\sqrt{-1}\kappa_{\nu} r \lambda/2} \cW_{\nu} p_{\mu}.
\end{eqnarray*}

\subsection{Reformulation of $\cW_{\mu}$ in terms of Schur polynomials}
\label{sec:W}

Recall
\begin{eqnarray}
&& \cW_{\mu} = q^{\kappa_{\mu}/4}\prod_{1 \leq i < j \leq l(\mu)}
\frac{[\mu_i - \mu_j + j - i]}{[j-i]}
\prod_{i=1}^{l(\mu)} \prod_{v=1}^{\mu_i} \frac{1}{[v-i+l(\mu)]},
\end{eqnarray}
where
\begin{align*}
q & = e^{\pi\sqrt{-1}\lambda}, & [m] & = q^{m/2} - q^{-m/2}.
\end{align*}
By the results in \cite{Zho1} we have
\begin{eqnarray}
\cW_{\mu}(q)
& = & q^{-|\mu|/2}s_{\mu}(1, q^{-1}, q^{-2}, \dots) \label{eqn:Wmu1} \\
& = & (-1)^{|\mu|} q^{\kappa_{\mu}/2+|\mu|/2}s_{\mu}(1, q, q^2, \dots). \label{eqn:Wmu2}
\end{eqnarray}

\subsection{One-partition Hodge integrals and KP hierarchy}

Now we can state our first main result:

\begin{theorem}
The generating series $G^{\bullet}(r; \lambda; p)$ is a $\tau$-function for the KP hierarchy in
variables $(x_1, x_2, \dots) = (p_1, \frac{1}{2}p_2, \dots)$.
\end{theorem}

\begin{proof}
By (\ref{eqn:Wmu1}) and (\ref{eqn:Wmu2}) we have
\begin{eqnarray*}
R^{\bullet}(r; \lambda; p)
& = & \sum_{\nu} e^{\kappa_{\nu}r/2}
s_{\nu}(x)  s_{\nu}(q^{-1/2}, q^{-3/2}, \dots) \\
& = & \sum_{\nu} e^{\kappa_{\nu}(r+1)/2}
s_{\nu}(x)  s_{\nu}(-q^{1/2}, -q^{3/2}, \dots).
\end{eqnarray*}
Now by (\ref{eqn:Skew}) and (\ref{eqn:K}),
we have
\begin{eqnarray*}
&& q^{\kappa_{\nu}(\tau+1)/2} s_{\nu}(x)
= \langle 0| Y_+(x_1, x_2, \dots)q^{(r+1)K}|\nu\rangle, \\
&&  s_{\nu}(-q^{1/2}, -q^{3/2}, \dots)
= \langle\nu|Y_-(-q^{1/2}, -q^{3/2}, \dots)|0\rangle,
\end{eqnarray*}
Therefore,
\begin{eqnarray*}
&& R^{\bullet}(r; \lambda; p) \\
& = & \sum_{\nu}
\langle 0| Y_+(x_1, x_2, \dots)q^{(r+1)K}|\nu\rangle
\cdot \langle\nu|Y_-(-q^{1/2}, -q^{3/2}, \dots)|0\rangle \\
& = & \langle 0| Y_+(x_1, x_2, \dots)q^{(r+1)K}Y_-(-q^{1/2}, -q^{3/2}, \dots)|0\rangle.
\end{eqnarray*}
It follows by Proposition \ref{prop:KP} that $R^{\bullet}(r; \lambda; p)$ is
a $\tau$-function of the KP hierarchy in variables $(p_1, \frac{1}{2}p_2, \dots)$,
and so is $G^{\bullet}(r; \lambda; p)$ by (\ref{eqn:MV}).
\end{proof}

\section{Two-Partition Hodge Integrals and $2$-Toda Hierarchy}

\label{sec:Two}

\subsection{Two-partition Hodge Integrals}

For a pair of partitions $(\mu^+, \mu^-) \in \cP^2_+$ (one of which might be empty),
consider the following generating series of Hodge integrals:
\begin{eqnarray*}
&& G_{\mu^+, \mu^-}(x, y) \\
& = & - \frac{\sqrt{-1}^{l(\mu^+)+l(\mu^-)}}{z_{\mu^+} \cdot z_{\mu^-}} \\
&& \cdot \sum_{g \geq 0} \lambda^{2g-2} \int_{\Mbar_{g, l(\mu^+)+l(\mu^-)}}
\frac{\Lambda_{g}^{\vee}(x)\Lambda^{\vee}_{g}(y)\Lambda_{g}^{\vee}(-x - y)}
{\prod_{i=1}^{l(\mu^+)} \frac{x}{\mu_i^+} \left(\frac{x}{\mu^+_i} - \psi_i\right)
\prod_{j=1}^{l(\mu^-)} \frac{y}{\mu_i^-}\left(\frac{y}{\mu^-_j} - \psi_{j+l(\mu^+)}\right)} \\
&& \cdot \left[xy(x+y)\right]^{l(\mu^+)+l(\mu^-)-1}
\cdot \prod_{i=1}^{l(\mu^+)} \frac{\prod_{a=1}^{\mu^+_i-1}
\left( \mu^+_iy + a x\right)}{\mu_i^+! x^{\mu_i^+-1}}
\cdot \prod_{i=1}^{l(\mu^-)} \frac{\prod_{a=1}^{\mu^-_i-1}
\left( \mu_i^- x + a y\right)}{\mu_i^-! y^{\mu_i^--1}}.
\end{eqnarray*}
The following formula was conjectured by the author \cite{Zho2} and proved in joint work with
C.-C. Liu and K. Liu \cite{LLZ3}:
\begin{eqnarray} \label{eqn:Zhou}
&& G^{\bullet}(r; p+,p^-) = R^{\bullet}(r; p^+, p^-),
\end{eqnarray}
where
\begin{eqnarray*}
G^{\bullet}(r; p^+, p^-) & = &
\exp \left(
\sum_{(\mu^+, \mu^-) \in \cP_+^2} G_{\mu^+, \mu^-}(r)p^+_{\mu^+}p^-_{\mu^-}\right), \\
R^{\bullet}(r; p^+, p^-) & = & \sum_{\mu^{\pm},\nu^{\pm}}
\frac{\chi_{\nu^+}(\mu^+)}{z_{\mu^+}} \frac{\chi_{\nu^-}(\mu^-)}{z_{\mu^-}}
e^{\sqrt{-1}(\kappa_{\nu^+} r  + \frac{\kappa_{\nu^-}}{r})\lambda/2}
\cW_{\nu^+, \nu^-} p^+_{\mu^+}p^-_{\mu^-}.
\end{eqnarray*}

\subsection{Reformulation of $\cW_{\mu, \nu}$ in terms of skew Schur functions}
\label{sec:W2}

Recall
\begin{eqnarray}
&& \cW_{\mu, \nu} = q^{|\nu|/2} \cW_{\mu} \cdot s_{\nu}(\cE_{\mu}(q,t)),
\end{eqnarray}
where\begin{eqnarray}
&& \cE_{\mu}(q,t) = \prod_{j=1}^{l(\mu)} \frac{1+q^{\mu_j-j}t}{1+q^{-j}t}
\cdot \left(1 + \sum_{n=1}^{\infty}
\frac{t^n}{\prod_{i=1}^n (q^i-1)}\right).
\end{eqnarray}
We have proved in \cite{Zho2}:
\begin{eqnarray} \label{eqn:Keynu}
&& s_{\nu}(\cE_{\mu}(q,t))
= (-1)^{|\nu|} q^{\kappa_{\nu}/2} \sum_{\rho} q^{-|\rho|}
\frac{s_{\mu/\rho}(1, q, q^2, \dots)}{s_{\mu}(1, q, q^2, \dots)}
s_{\nu/\rho}(1, q, q^2, \dots).
\end{eqnarray}
and
\begin{eqnarray} \label{eqn:Key}
&& \cW_{\mu, \nu}(q)
= (-1)^{|\mu|+|\nu|}
q^{\frac{\kappa_{\mu}+\kappa_{\nu}+|\mu|+|\nu|}{2}}
\sum_{\rho} q^{-|\rho|} s_{\mu/\rho}(1, q, \dots)s_{\nu/\rho}(1, q, \dots).
\end{eqnarray}

\subsection{Two-partition Hodge integrals and $2$-Toda hierarchy}

\begin{theorem}
Let
$$\tau_n = q^{(r+\frac{1}{r}+2)\frac{n(4n^2-1)}{24}}
G^{\bullet}(r; \lambda; p^+(q^{(\frac{1}{r}+1)n}x^+), p^-(q^{(\frac{1}{r}+1)n}x^-)).$$
Then $\{\tau_n\}$ is a sequence of $\tau$-functions of the $2$-Toda hierarchy
in two sequences of variables
$(t_1^\pm, t_2^{\pm}, \dots) = (p_1(x^{\pm}), \frac{1}{2}p_2(x^{\pm}), \dots)$.
\end{theorem}

\begin{proof}
By (\ref{eqn:Key}) we have
\begin{eqnarray*}
&& R^{\bullet}(r; p(x^+), p(x^-)) \\
& = & \sum_{|\mu^{\pm}|=|\nu^{\pm}|}
\frac{\chi_{\nu^+}(\mu^+)}{z_{\mu^+}}
\frac{\chi_{\nu^-}(\mu^-)}{z_{\mu^-}}
e^{\sqrt{-1}(\kappa_{\nu^+}r + \kappa_{\nu^-}/r)\lambda/2} \cW_{\nu^+, \nu^-}(q)
p^+_{\mu^+}p^-_{\mu^-} \\
& = & \sum_{\nu^{\pm}} s_{\nu^+}(x^+) q^{\kappa_{\nu^+}r/2}
\cW_{\nu^+, \nu^-}(q) q^{\kappa_{\nu^-}/(2r)}s_{\nu^-}(x^-) \\
& = & \sum_{\nu^{\pm}} s_{\nu^+}(x^+) s_{\nu^-}(x^-)
e^{\kappa_{\nu^+}(r+1)/2+\kappa_{\nu^-}(\frac{1}{r}+1)/2} \\
&& \cdot
 \sum_{\rho} s_{\nu^+/\rho}(-q^{1/2}, -q^{3/2}, \dots)s_{\nu^-/\rho}(-q^{1/2}, -q^{3/2}, \dots)\\
& = &  \sum_{\rho} \sum_{\nu^+} e^{\kappa_{\nu^+}(r+1)/2}
 s_{\nu^+}(x^+)  s_{\nu^+/\rho}(-q^{1/2}, -q^{3/2}, \dots) \\
&& \cdot \sum_{\nu^-} e^{\kappa_{\nu^-}(\frac{1}{r}+1)/2}
s_{\nu^-}(x^-)s_{\nu^-/\rho}(-q^{1/2}, -q^{3/2}, \dots).
\end{eqnarray*}
Now by (\ref{eqn:Skew}) and (\ref{eqn:K}) we have
\begin{eqnarray*}
&& q^{\kappa_{\nu^+}(r+1)/2} s_{\nu^+}(x^+)
= \langle 0| Y_+(x_1^+, x_2^+, \dots)q^{(r+1)K}|\nu^+\rangle, \\
&&  s_{\nu^+/\rho}(-q^{1/2}, -q^{3/2}, \dots)
= \langle\nu^+|Y_-(-q^{1/2}, -q^{3/2}, \dots)|\rho\rangle, \\
&& q^{\kappa_{\nu^-}(\frac{1}{r}+1)/2} s_{\nu^-}(x^-)
= \langle \nu^-|q^{(\frac{1}{r}+1)K}Y_-(x_1^-, x_2^-, \dots)|0\rangle, \\
&& s_{\nu^-/\rho}(-q^{1/2}, -q^{3/2}, \dots)
= \langle \rho|Y_+(-q^{1/2}, -q^{3/2}, \dots)|\nu^-\rangle,
\end{eqnarray*}
therefore
\begin{eqnarray*}
&& R^{\bullet}(r; p(x^+), p(x^-)) \\
& = & \sum_{\rho, \nu^+, \nu^-}
\langle 0| Y_+(x_1^+, x_2^+, \dots)q^{(r+1)F_2}|\nu^+\rangle
\cdot \langle\nu^+| Y_-(-q^{1/2}, -q^{3/2}, \dots)|\rho\rangle \\
&& \cdot \langle \rho| Y_+(-q^{1/2}, -q^{3/2}, \dots)|\nu^-\rangle
\cdot \langle \nu^-|q^{(\frac{1}{r}+1)K} Y_-(x_1^-, x_2^-, \dots)|0\rangle \\
& = & \langle 0| Y_+(x_1^+, x_2^+, \dots)q^{(r+1)K}Y_-(-q^{1/2}, -q^{3/2}, \dots)\\
&& \cdot Y_+(-q^{1/2}, -q^{3/2}, \dots)q^{(\frac{1}{r}+1)K} Y_-(x_1^-, x_2^-, \dots)|0\rangle.
\end{eqnarray*}
Hence by Proposition \ref{prop:Toda},
$R^{\bullet}(\tau; p(x^+), p(x^-))$ and so $G^{\bullet}(\tau; p(x^+), p(x^-))$ by (\ref{eqn:Zhou})
is the $\tau_1$ of a sequence $\{\tau_n: n \geq 1\}$
of $\tau$-functions of the $2$-Toda hierarchy
in the variables $(t_1^\pm, t_2^{\pm}, \dots) = (p_1(x^{\pm}), \frac{1}{2}p_2(x^{\pm}), \dots)$.
By (\ref{eqn:Translation}), (\ref{eqn:H}), and (\ref{eqn:Comm}),
we have
\begin{eqnarray*}
\tau_n & = &
\langle n| Y_+(x^+)q^{(r+1)K}Y_-(-q^{1/2}, -q^{3/2}, \dots)\\
&& \cdot Y_+(-q^{1/2}, -q^{3/2}, \dots)q^{(\frac{1}{r}+1)K} Y_-(x^-)|n\rangle \\
& = & \langle 0| R^{-n}Y_+(x^+)q^{(r+1)K}Y_-(-q^{1/2}, -q^{3/2}, \dots)\\
&& \cdot Y_+(-q^{1/2}, -q^{3/2}, \dots)q^{(\frac{1}{r}+1)K} Y_-(x^-)R^n |0\rangle \\
& = & \langle 0| Y_+(x^+)R^{-n}q^{(r+1)K}R^n Y_-(-q^{1/2}, -q^{3/2}, \dots)\\
&& \cdot Y_+(-q^{1/2}, -q^{3/2}, \dots)R^{-n}q^{(\frac{1}{r}+1)K}R^n  Y_-(x^-)|0\rangle \\
& = & \langle 0| Y_+(x^+)q^{(r+1)(K+nH+\frac{n^2}{2}\alpha_0+\frac{n(4n^2-1)}{24})}
Y_-(-q^{1/2}, -q^{3/2}, \dots)\\
&& \cdot Y_+(-q^{1/2}, -q^{3/2}, \dots)
q^{(\frac{1}{r}+1)(K+nH+\frac{n^2}{2}\alpha_0+\frac{n(4n^2-1)}{24})}
Y_-(x^-)|0\rangle \\
& = & q^{(r+\frac{1}{r}+2)\frac{n(4n^2-1)}{24}}
\langle 0| Y_+(q^{(\frac{1}{r}+1)n}x^+)q^{(r+1)K}
Y_-(-q^{1/2}, -q^{3/2}, \dots)\\
&& \cdot Y_+(-q^{1/2}, -q^{3/2}, \dots)
q^{(\frac{1}{r}+1)K}
Y_-(q^{(\frac{1}{r}+1)n}x^-)|0\rangle \\
& = & q^{(r+\frac{1}{r}+2)\frac{n(4n^2-1)}{24}}
G^{\bullet}(r; \lambda; p^+(q^{(\frac{1}{r}+1)n}x^+), p^-(q^{(\frac{1}{r}+1)n}x^-)).
\end{eqnarray*}
This completes the proof.
\end{proof}

\begin{remark}
From the above proof it is easy to see that
\begin{eqnarray}
&& \cW_{\nu^+, \nu^-} = \langle \nu^+|q^KY_+(-q^{\frac{1}{2}}, - q^{\frac{3}{2}} , \dots)
Y_-(-q^{\frac{1}{2}}, - q^{\frac{3}{2}} , \dots) q^K |\nu^- \rangle.
\end{eqnarray}
\end{remark}

\section{Relative Invariants and Integrable Hierarchies}

\label{sec:Relative}

In the above we have established the relationships between some special Hodge integrals
with some well-known integrable hierarchies.
We are partly inspired by remarks made by Dijkgraaf in a recent talk on open string theory:
The disc gives KP hierarchy,
the cylinder gives $2$-Toda equation.
(The author thanks Chiu-Chu Liu for communicating Dijkgraaf's remarks to him.)
While open string invariants do not in general have rigorous mathematical formulations,
the relative string invariants can be rigorously defined using the moduli spaces of J. Li \cite{Li}.
As the discussions in \cite{Li-Son} indicates,
there might be a duality between the open and relative invariants.
In this section we will deal with some relative invariants which has been computed as the byproducts
of the proof of (\ref{eqn:MV}) and (\ref{eqn:Zhou}) in \cite{LLZ1} and \cite{LLZ3} respectively.

\subsection{Relative invariants of the resolved conifold and KP hierarchy}

The resolved conifold is the space $\cO(-1) \oplus \cO(-1) \to \bP^1$.
Consider the $\bC^*$-action
$$t \cdot [z^0:z^1] = [tz^0: z^1]$$
on $\bP^1$. It has two fixed points
$p_0= [0:1]$ and $p_1=[1:0]$.
For a partition $\mu$ of $d>0$
let $\Mbar_{g, 0}(\bP^1, \mu)$ be the moduli
space of morphisms relative to $p_1$ with ramification type $\mu$.
It is a separated, proper Deligne-Mumford stack with
a perfect obstruction theory of virtual dimension
$$r_{g, \mu}=2g-2 + |\mu| + l(\mu),$$
so it has a virtual fundamental class of degree $r_{g, \mu}$.

One can define two bundles $V_D$ and $V_{D_d}$ on $\Mbar_{g, 0}(\bP^1, \mu)$ as follows.
Let
$$
\pi:\mathcal{U}_{g,\mu}\to\Mbar_{g,0}(\bP^1,\mu)
$$
be the universal domain curve, and let
$$
P:\mathcal{T}_{g,\mu}\to\Mbar_{g,0}(\bP^1,\mu)
$$
be the universal target.
There is an evaluation map
$$
F:\mathcal{U}_{g,\mu}\to \mathcal{T}_{g,\mu}
$$
and a contraction map
$$
\tilde{\pi}:\mathcal{T}\to \bP^1.
$$
Let $\mathcal{D}_{g,\mu}\subset \mathcal{U}_{g,\mu}$ be the
divisor corresponding to the $l(\mu)$ marked points.
Define
\begin{eqnarray*}
V_D&=&R^1\pi_*(\cO_{\mathcal{U}_{g,\mu}}(-\mathcal{D}_{g,\mu}) )\\
V_{D_d}&=&R^1\pi_* \tilde{F}^*\cO_{\bP^1}(-1),
\end{eqnarray*}
where $\tilde{F}=\tilde\pi\circ F:\mathcal{U}_{g,\mu}\to \bP^1$.
The bundle
$$
V_{g, \mu} =V_D\oplus V_{D_d}
$$
is a vector bundle of rank $r_{g, \mu}=2g-2+d+l(\mu)$.
(For more details, see \cite{LLZ1, LLZ2}.)
Define
\begin{eqnarray*}
&& K_{g, \mu} = \int_{\Mbar_{g,0}(\bP^1,\mu)} e(V_{g, \mu}), \\
&& K_{\mu}(\lambda) = \sum_{g \geq 0} \lambda^{2g-2} K_{g, \mu}, \\
&& K(\lambda; p) = \sum_{|\mu| > 0} K_{\mu}(\lambda)p_{mu}, \\
&& K^{\bullet}(\lambda;p) = e^{K(\lambda; p)}.
\end{eqnarray*}

The $\bC^*$-action induces induces $\bC^*$-actions on
$\Mbar_{g, 0}(\bP^1, \mu)$ and $V$,
and one can compute $K_{g, \mu}$ by localization method \cite{Gra-Pan}.
This has been exploited in \cite{LLZ1, LLZ2} to prove (\ref{eqn:MV}).
By the localization calculations there with $r = -1$
($\tau = -1$ in the notation there),
one has
\begin{eqnarray*}
K_{\mu}(\lambda) & = & G_{\mu}(-1; \lambda),
\end{eqnarray*}
hence by (\ref{eqn:MV})
\begin{eqnarray*}
K^{\bullet}(\lambda;p) & = & G^{\bullet}(-1;\lambda;p)
= R^{\bullet}(-1;\lambda; p) \\
& = & \langle 0| Y_+(x_1, x_2, \dots)Y_-(-q^{1/2}, -q^{3/2}, \dots)|0\rangle.
\end{eqnarray*}
Therefore,
$K^{\bullet}(\lambda;p)$ is a $\tau$-function of the KP hierarchy.

\begin{eqnarray*}
&& K_{g, \mu} = \int_{\Mbar_{g,0}(\bP^1,\mu)} e(V), \\
&& K_{\mu}(\lambda) = \sum_{g \geq 0} \lambda^{2g-2} K_{g, \mu}, \\
&& K(\lambda; p) = \sum_{|\mu| > 0} K_{\mu}(\lambda)p_{mu}, \\
&& K^{\bullet}(\lambda;p) = e^{K(\lambda; p)}.
\end{eqnarray*}

The $\bC^*$-action induces induces $\bC^*$-actions on
$\Mbar_{g, 0}(\bP^1, \mu)$ and $V_{g, \mu}$,
and one can compute $K_{g, \mu}$ by localization method \cite{Gra-Pan}.
This has been exploited in \cite{LLZ1, LLZ2} to prove (\ref{eqn:MV}).
By the localization calculations there with $r = -1$
($\tau = -1$ in the notation there),
one has
\begin{eqnarray*}
K_{\mu}(\lambda) & = & G_{\mu}(-1; \lambda),
\end{eqnarray*}
hence by (\ref{eqn:MV})
\begin{eqnarray*}
K^{\bullet}(\lambda;p) & = & G^{\bullet}(-1;\lambda;p)
= R^{\bullet}(-1;\lambda; p) \\
& = & \langle 0| Y_+(x_1, x_2, \dots)Y_-(-q^{1/2}, -q^{3/2}, \dots)|0\rangle.
\end{eqnarray*}
Therefore,
$K^{\bullet}(\lambda;p)$ is a $\tau$-function of the KP hierarchy.

\subsection{Relative invariants of a toric Fano surface and the $2$-Toda hierarchy}

Let $X$ be the toric surface obtained by blowing up $(p_1, p_1)$ on
$\bP^1 \times \bP^1$.
Denote by $D^+$ and $D^-$ the strict transform of the divisors $\bP^1 \times \{p_1\}$ and
$\{p_1\} \times \bP^1$ respectively.
For $(\mu^+,\mu^-)\in\cP^2_+$,
let $\Mbar_{g,0}(X; \mu^+,\mu^-)$ be the moduli spaces of morphisms which have ramification
type $\mu^\pm$ along $D^\pm$.
It is a separated, proper Deligne-Mumford stack with
a perfect obstruction theory of virtual dimension
$$
r_{g, \mu^+, \mu^-}=g-1 +|\mu^+|+l(\mu^+)+|\mu^-|+l(\mu^-),
$$
so it has a virtual fundamental class of degree $r_{g, \mu^+, \mu^-}$.

Let
$$
\pi:\gmu
{U}\to\MgX
$$
be the universal domain curve, and let
$$
  P:\gmu{T}\to\MgX
$$
be the universal target.
There is an evaluation map
$$
F:\gmu{U}\to \gmu{T}
$$
and a contraction map
$$
\tilde{\pi}: \gmu{T}\to X.
$$
Let $\gmu{D}\subset \gmu{U}$ be the
divisor corresponding to the $l(\mu^+)+l(\mu^-)$ marked points.
Define
$$
\gmu{V}=R^1\pi_*\left( \tilde{F}^*\cO_X(-D_1-D_3)\otimes
\cO_{\gmu{U}}(-\gmu{D})\right)
$$
where $\tilde{F}=\tilde{\pi}\circ F:\gmu{U}\to X$.
Now $V_{g, \mu^+, \mu^-}\to \MgX$ is a vector bundle of rank
$$
r_{g, \mu^+, \mu^-} =g-1+|\mu^+|+l(\mu^+)+|\mu^-|+l(\mu^-).
$$
(For details see \cite{LLZ3}.)
Define
\begin{eqnarray*}
&& K_{g, \mu^+, \mu^-} = \int_{\Mbar_{g,0}(X;\mu^+, \mu^-)} e(V_{g, \mu^+, \mu^-}), \\
&& K_{\mu^+, \mu^-}(\lambda) = \sum_{g \geq 0} \lambda^{2g-2} K_{g, \mu^+, \mu^-}, \\
&& K(\lambda; p(x^+), p(x^-)) = \sum_{(\mu^+, \mu^-) \in \cP^2_+} K_{\mu^+, \mu^-}(\lambda)p_{\mu^+}(x^+)p_{\mu^-}(x^-), \\
&& K^{\bullet}(\lambda;p(x^+), p(x^-)) = e^{K(\lambda; p(x^+), p(x^-))}.
\end{eqnarray*}

The $(\bC^*)^2$-action induces $(\bC^*)^2$-actions on
$\Mbar_{g, 0}(X, \mu^+, \mu^-)$ and $V_{g, \mu^+, \mu^-}$,
and one can compute $K_{g, \mu^+, \mu^-}$ by localization method \cite{Gra-Pan}.
This has been exploited in \cite{LLZ1, LLZ2} to prove (\ref{eqn:Zhou}).
As a byproduct the following identity has been proved in \cite{LLZ3}:
\begin{eqnarray*}
&& K^\bullet_{\mu^+, \mu^-}(\lambda)
= \sum_{\nu^\pm}
\frac{\chi_{\nu^+}(C_{\mu^+})}{z_{\mu^+}}
\cdot
\cW_{\nu^+, \nu^-}
\cdot  \frac{\chi_{\nu^-}(C_{\mu^-})}{z_{\mu^-}},
\end{eqnarray*}
hence
\begin{eqnarray*}
&& K^{\bullet}(\lambda;p(x^+), p(x^-)) \\
& = & \sum_{\nu^\pm}s_{\nu^+}(x^+) \cW_{\nu^+, \nu^-} s_{\mu^-}(x^-) \\
& = & \langle 0| Y_+(x^+)q^KY_-(-q^{1/2}, -q^{3/2}, \dots)
Y_+(-q^{1/2}, -q^{3/2}, \dots)q^KY_-(x^-)|0\rangle.
\end{eqnarray*}
In this case,
\begin{eqnarray*}
&& \langle n| Y_+(x^+)Y_-(-q^{1/2}, -q^{3/2}, \dots)
Y_+(-q^{1/2}, -q^{3/2}, \dots)Y_-(x^-)|n\rangle \\
& = & \langle 0| R^{-n}Y_+(x^+)q^KY_-(-q^{1/2}, -q^{3/2}, \dots)
Y_+(-q^{1/2}, -q^{3/2}, \dots)q^KY_-(x^-)R^n|0\rangle \\
& = & q^{\frac{n(4n^2-1}{12}}\langle 0| Y_+(q^nx^+)Y_-(-q^{1/2}, -q^{3/2}, \dots)
Y_+(-q^{1/2}, -q^{3/2}, \dots)Y_-(q^nx^-)|0\rangle\\
& = & q^{\frac{n(4n^2-1}{12}} K^{\bullet}(\lambda;p(q^nx^+), p(q^nx^-)).
\end{eqnarray*}
Therefore,
$\tau_n = q^{\frac{n(4n^2-1}{12}}K^{\bullet}(\lambda;p(q^nx^+), p(q^nx^-))$ is a sequence of $\tau$-functions of the $2$-Toda hierarchy.

\subsection{A conjecture on relative invariants}
Based on the above examples,
we formulate the following

\begin{conjecture}
Let $X$ be a toric Fano surface.
Then the generating series of suitably defined invariants on moduli spaces
relative to a toric invariant divisor $D$ is a $\tau$-function of the KP hierarchy,
and the the generating series of suitably defined invariants on moduli spaces
relative to two disjoint toric invariant divisors $D^+$ and $D^-$ gives rise to
a sequence of $\tau$-functions of the $2$-Toda hierarchy.
\end{conjecture}

{\bf Acknowledgement}.
The author thanks Professors Chiu-Chu Melissa Liu and Kefeng Liu for collaborations
which lead to this work.
The research in this work is partly supported by research grants from NSFC and Tsinghua University.

\end{document}